\documentclass{article}
\usepackage{dirtytalk}
\usepackage[margin=3cm]{geometry}
\usepackage[utf8]{inputenc}
\usepackage[english]{babel}
\usepackage{amsthm}
\usepackage{amsmath}
\usepackage{amssymb}
\usepackage{amsfonts}
\usepackage{xspace}
\usepackage{hyperref}
\usepackage{mathtools}
\usepackage{tikz-cd}
\usepackage{hyperref}
\usepackage[style=alphabetic]{biblatex}
\hypersetup{
    colorlinks=true,
    linkcolor=blue,
    filecolor=magenta,      
    urlcolor=cyan,
}
\usepackage[title]{appendix}
\theoremstyle{definition}
\newtheorem{definition}{Definition}[section]
\theoremstyle{plain}
\newtheorem{theorem}{Theorem}[section]
\newtheorem{lemma}[theorem]{Lemma}
\newtheorem{corollary}[theorem]{Corollary}
\theoremstyle{remark}
\newtheorem{remark}[theorem]{Remark}
\newtheorem{example}[theorem]{Example}

\newcommand{\ie}{\textit{i.e.}\xspace}

\newcommand{\N}{\mathbb{N}}

\title{A dynamical version of Silverman-Tate's height inequality}
\author{Debam Biswas\ \&\ Zhelun Chen}
\date{}

\begin{document}

\maketitle
\section*{Contact}
Address(first author):\ \emph{Debam Biswas,\ Department of Mathematics,\ University of Regensburg\\
Universitatstrasse 31, 93053, Regensburg}\\
Address(second author):\ \emph{Zhelun Chen,\ Mathematical Institute,\ Leiden University,\\
Einsteinweg 55,\ 2333 CC Leiden,\ The Netherlands}\\
Email(first author):\ \href{mailto:debambiswas@gmail.com}{debambiswas@gmail.com}\\
Email(second author):\ \href{mailto:zl.chen1729@gmail.com}{zl.chen1729@gmail.com}

\begin{abstract}
In the paper "Uniformity of Mordell-Lang" by Vesselin Dimitrov, Philipp Habegger and Ziyang Gao (\cite{dimitrov2020uniformity}, they use Silverman's Height Inequality and they give a proof of the same which makes use of Cartier divisors and hence drops the flatness assumption of structure morphisms of compactified abelian schemes. However, their proof makes use of Hironaka's theorem on resolution of singularities which is unknown for fields of positive characteristic. We try to slightly modify their ideas, use blow-ups in place of Hironaka's theorem to make the proof effective for any fields with product formula where heights can be defined and any normal quasi.projective variety as a base. We work in the more general set-up of dynamical systems. As an application we prove certain variants of the specialisation theorems in \cite{Sil83} with restricted hypotheses in higher dimensional bases.
\end{abstract}
\section*{Keywords:}
\textbf{2020 Mathematics Subject Classification :\ 11G50,\ 14G40}

\section{Introduction}
\subsection{A short history}
Let $\mathcal{A}$ be a family of abelian varieties parametrized by a normal, irreducible and quasi-projective variety $S$ over a field $k$ with Product Formula. Consider the compactified familily $\overline{\mathcal{A}}\to \overline{S}$. There are three kinds of height
functions associated to such a family:  the  Weil heights $h_{\overline{\mathcal{A}}}$ and $h_{\overline{S}}$ on the family $\overline{\mathcal{A}}$ and on the  base $\overline{S}$ respectively, and the N\'eron-Tate
height $h_{s}$ on each fiber of the family $\mathcal{A}\to S$. It is natural to ask how these heights are related.

In 1983, Silverman gave an answer to this question: he established and proved an inequality \cite[Theorem A]{Sil83}, which will be called the Height Inequality in this paper, in the situation where the family $\overline{\mathcal{A}}\to \overline{S}$ is flat and $\overline{S}$ is smooth over a field satisfying the Product formula.\ In the same year Tate improved a consequence of the Height Inequality, namely \cite[Theorem B]{Sil83} concerning the limit behaviour of the heights involved in the situation where $\overline{\mathcal{A}}\to \overline{S}$ is an elliptic surface (\cite[\S 7]{Sil83}) and the parameter space is a curve.

When the base space is a smooth curve,\ the compactified structure morphism is always dominant and hence flat.\ However the same is not true for parameter spaces of arbitary dimension.\ In the recent paper \cite[Appendix A]{dimitrov2020uniformity} the authors proved the Height Inequality without flatness assumption on the family $\overline{\mathcal{A}}\to \overline{S}$ but $S$ is assumed to be regular. They use it (combined with other more sophisticated methods) to prove an old conjecture of Mazur on uniformity in the Mordell-Lang conjecture.\ In their proof without flatness,\ the authors use Hironaka's famous theorem on resolution of singularities.\ As Hironaka's theorem is unknown in fields of characteristic 0,\ this makes the proof in \cite{dimitrov2020uniformity} only true for fields of characteristic 0.

In this article we adapt the above scenario in the more general framework of \emph{polarised dynamical systems}, one of the principle examples of which are families of abelian varieties together with an even ample line bundle.  We replace the arguments involving Hironaka's theorem with simple blow-up techniques and slightly modifying the arguments,\ we obtain a different proof. So in comparison to the previously obtained results in \cite[Theorem A]{Sil83} and \cite[Theorem 3.1]{SilCall}, we are able to relax the smoothness assumption on the base space to only normality over fields of arbitary characteristic satisfying a Product formula and relax the flatness assumption of the structure morphism to dominancy. In the last section we treat the problem of scarcity of exceptional points and obtain a result similar to \cite[Theorem 4.1]{SilCall} allowing our base space to be projective space of any dimension but restricting the variation of the base parameter only on closed sub-varieties.\\
Lastly we will like to talk about applications of the Silverman-Tate height inequality.\ The inequality was used originally in \cite{Sil83} to compare generic Ne\'ron-Tate heights to the fiberwise ones via a limiting process. This relation was in turn used to show a specialisation theorem as in Theorem C of \cite{Sil83}. We prove similar partial results in the case where the base-space is a higher dimensional projective space and not only one dimensional curves.  As mentioned earlier,\ this inequality is used together with other sophisticated methods to prove Uniformity in Mordell-Lang in the recent paper \cite{dimitrov2020uniformity}. Another application of the Silverman Tate inequality is in the proof of the Geometric Bogomolov Conjecture for the function field of a curve defined over $\overline{\mathbb{Q}}$ in the paper \cite{gao2019heights} by Gao and Habegger.
\subsection{Dynamical systems and the main theorem}
\noindent\textbf{Convention}. 
\begin{itemize}
\item By a \textbf{variety} over $k$ we mean an integral scheme which is of finite type over a field~$k$.
\item A morphism $X\to S$  of scheme is called \textbf{projective} if there is a closed immersion $X\xhookrightarrow{} \mathbb{P}^n_S$ over~$S$ for some  $n\geq 0$. (See \cite[Summary 13.71]{gortz2010algebraic} for the comparison of different conventions.)
\end{itemize} 

Following \cite[\S 6.1]{YZ21}, we introduce
\begin{definition}
A \textbf{polarized dynamical system} is a tuple $(X,f,L,q)_{/S}$ consisting of the following data:
\begin{enumerate}
\item There is a morphism $\pi: X\to S$ of integral schemes which is flat and projective.
\item $f: X\to X$ is an $S$-morphism.
\item $L$ is a line bundle on $X$, relatively ample with respective to $\pi$, such that
\[
f^*L\cong L^{\otimes q}\quad \text{(polarization)}
\] for some natural number $q>1$.
\end{enumerate}
\end{definition}
For the basics of relative ample line bundles we refer to \cite[\S 1.7]{Laz04}. We recall that relative ampleness can be characterized fiberwise (when $S$ is a $k$-variety). From the polarization we see that $f$ is a finite morphism of degree $q$.

From now on  $k$ denotes a \textbf{field with product formula} of any characteristics.  For the sake of concreteness one can think of $k$ to be a number field or a function field of a normal curve (1-dimensional variety) over an algebraically closed field. If $X$ is a projective variety over $k$ and $L$ is a line bundle on $X$, then one can define a (real-valued) height function $h_L$ associated to $L$ on the set of algebraic points $X(\overline{k})$ via the so-called Weil's height machine. We refer to \cite[\S 3 \& 4]{lang2013fundamentals} for the details. We also write $h_X$ for $h_{\mathcal{O}(1)}$, where $\mathcal{O}(1)$ is the ample line bundle on $X$ associated to the projective embedding $X\xhookrightarrow{} \mathbb{P}^n_k$ for some $n\geq 0$.  

Let us introduce the heights in the relative setting; this is achieved fiberwise as follows.  Consider a projective morphism $\pi: X\to S$ of integral schemes, fix a projective compactification $S\subset \overline{S}\subset \mathbb{P}^m_k$ (so that $S$ is an open subscheme of $\overline{S}$ and the latter is a closed subscheme of $\mathbb{P}^m_k$, all defined over $k$), where $S$ is taken to be a normal quasi-projective variety over $k$. 

Like in \cite[Appendix A]{dimitrov2020uniformity}, we have a closed immersion 
$$X\xhookrightarrow{}S\times \mathbb{P}_k^n.$$ Together we have an immersion 
$$X\xhookrightarrow{}\mathbb{P}_k^m\times \mathbb{P}_k^n$$
and let us call by $\overline{X}$ the Zariski closure of $X$ under this immersion in the multi-projective space $\mathbb{P}_k^{m,n}:=\mathbb{P}_k^m\times\mathbb{P}_k^n$.

We have a very ample line bundle $\mathcal{O}(1,1)$ on the multiprojective space~$\mathbb{P}_k^{m,n}$.\ Let
\[\mathcal{L}:=\mathcal{O}(1,1)\mid_{\overline{X}}\ \text{and}\ L\cong \mathcal{L}\mid_X\]
A point in $\overline{X}(\overline{k})$ can be represented by a pair $(s,P)\in  \overline{S}(\overline{k})\times \mathbb{P}^n_k(\overline{k})$. Thus the height function $h_{\mathcal{L}}$  on $\overline{X}(\overline{k})\subset \mathbb{P}^{m,n}_k(\overline{k})$ can be described (up to bounded functions) by
$$h_{\mathcal{L}}(s,P)=h_{\overline{S}}(s)+h_n(P),$$ where $h_n$ denotes the standard Weil height  on $\mathbb{P}_k^n$.\ For notational convenience we also denote the restriction of $h_{\mathcal{L}}(\cdot)$ to $X$ as $h_L(\cdot)$. By choice, the line bundle $L$ is  relatively ample with respect to $\pi$ . Consider an $S$-morphism $f: X\to X$. We assume that there is a natural number $q>1$ such that 
\[
f^*L\cong L^{\otimes q}.
\] In this way we obtain a polarized dynamical system $(X,f,L,q)_{/S}$.
\begin{example}\label{abelian scheme as dynamical system}
Let $S$ be a quasi-projective normal variety, $X$ be an abelian scheme over $S$ and $L$ an even line bundle which is ample with respect to the structural morphism $X\to S$. Then by the theorem of the cube we see that $(X,[2],L,4)_{/S}$ is a polarized dynamical system in the above sense. That $X\to S$ is projective follows from the normality of $S$ (a theorem of Raynaud), cf. \cite[Remark 3.1]{dimitrov2020uniformity}.
\end{example}
Our main result is the following height inequality.
\begin{theorem}\label{dynamical ST height inequality}
Let $(X,f,L,q)_{/S}$ be a polarized dynamical system constructed as above. Then there exists a constant $c>0$ such that  for all $P\in X(\overline{k})$, we have 
$$|h_{L}(f(P))-qh_{L}(P)|\le c\cdot\max\{1,h_{\overline{S}}(\pi(P))\}$$.
\end{theorem}
Our result then immediately gives as a corollary an inequality relating the canonical heights of a polarised dynamical system with the Naive height as in \cite[Theorem 3.1]{SilCall}. Indeed suppose $(X,S,f)$ be a polarized dynamical system with $f^*L\cong L^{\otimes q}$. Then for all $t\in S(\overline{k})$ if we denote $L_t:=L|_{X_t}=\mathcal{L}|_{X_t}$ to be the restriction of $\mathcal{L}$ to the fiber $X_t$ over $t$, then we clearly have $f^*L_t\cong L_t^{\otimes q}$. Then a standard Tate's limiting argument together with the Height machine allows us to define a \emph{canonical height} as follows
\[
\hat{h}_X(P)=\lim_{n\to \infty} \frac{h_{\overline{X}}(f^{(n)}(P))}{q^n}\quad \text{for all algebraic points } P\in X(\overline{k}),
\] where the (fiberwise) Weil height $h_{\overline{X}}$ is induced by the $\pi$-ample even line bundle $\mathcal{L}$ considered before and $f^{(n)}(P)$ is the $n$-fold application of $f$ to $P$. Then we can state our corollary:
\begin{corollary}[Silverman-Tate's height inequality]\label{ST's height ineq}
Let $\pi: X\to S$ be a dynamical system with $f^*L\cong L^{\otimes q}$. Then there is a constant $c>0$ such that for all $P\in X(\overline{k})$, we have
\[
\left|\hat{h}_X(P)-h_{\overline{X}}(P)\right|\leq c\cdot\max\left(1,h_{\overline{S}}(\pi(P))\right).
\]
\end{corollary}
\begin{proof}
(Cf. proof of \cite[Theorem A.1]{dimitrov2020uniformity}) Pick any $P\in X(\bar{k})$. Let $\ell\geq m\geq 0$ be integers. Using telescope-sum and triangle inequality, we find
\[
\left|\frac{h(f^{(l)}(P))}{q^l}-\frac{h(f^{(m)}(P))}{q^m}\right|\leq \sum_{t=m}^{l-1}\left|\frac{h(f^{(t+1)}(P))-qh(f^{(t)}(P))}{q^{t+1}}\right|.
\] Applying Theorem \ref{dynamical ST height inequality} to the summands, we find a constant $c>0$ such that the right hand side is bounded by
\[
c\lambda\sum_{q=m}^{\ell-1}\frac{1}{4^{q+1}}\leq c\lambda\frac{1}{4^m},\quad \lambda:=\max(1,h_{\overline{S}}(\pi(P))).
\] This means that $(\frac{h(f^{(l)}(P))}{4^l})_{l\geq 1}$ is a Cauchy sequence with limit $\hat{h}_X(P)$, simply by definition of the N\'eron-Tate height. On taking $m=0$ and considering $\ell\to\infty$ we obtain Corollary \ref{ST's height ineq}.
\end{proof}

\textbf{Acknowledgement}. We are indebted to Prof. Walter Gubler for providing us with the main idea and helping us to see the possibility of the generalisation to normal bases over arbitary global fields. We are also indebted to Prof. Joseph Silverman for pointing the connection with \cite{SilCall}. The first author is grateful to Prof. Qing Liu for answering  questions regarding a generalisation. The first author would also like to thank Prof. Yuri Bilu for introducing him to the very recent paper \cite{dimitrov2020uniformity} as his potential master’s thesis and Ananyo Kazi for numerous fruitful discussions and for knowing Hartshorne by heart. The second author thanks Prof. Robin de Jong for the encouragement for the generalization of Silverman-Tate's height inequality to the dynamical setting and for notifying the relation of Silverman's specialization theorem and the \emph{uniform torsion conjecture}, cf. \cite{Holmes}.
\section{Extension of the polarising morphism}

To prove Theorem \ref{dynamical ST height inequality}, we first need to handle the problem that the self-map $f: X\to X$ usually does not extend  to an ordinary $\overline{S}$-morphism  $\overline{X}\to \overline{X}$ of the compactified family  (but merely a rational map).  We bypass this problem by going to the products via the graph morphism and then identifying the image subscheme with $X$. Here we imitate \cite[Appendix A]{dimitrov2020uniformity}.

We consider the graph morphism associated to $f$,\ namely:\\
$$X\xhookrightarrow{\Gamma_{f}} X\times_S X\xhookrightarrow{}\overline{X}\times_{\overline{S}}\overline{X}.$$
As $X$ is separated over $S$,\ the first morphism is a closed immersion while the second one is an open immersion.\ We denote by $\overline{X}'$ the Zariski closure of $X$ under this immersion.\ We can identify $X\subset \overline{X}'$ as an open dense subscheme.\ Moreover we have the two canonical projections $p_1$ and $p_2$ on the product $\overline{X}\times_{\overline{S}}\overline{X}$ and note that by construction 
$$p_1\mid_{X}=\text{id}_{X}\ \text{and}\ p_2\mid_{X}=f.$$
To summarise,\ we have the following commutative diagram:\\
$$
\begin{tikzcd}
\overline{X}                & \overline{X}' \arrow[r, "{p_2}"] \arrow[l, "p_1"'] & \overline{X}      \\
 X \arrow[u, hook]  & X \arrow[u, hook] \arrow[swap, r,"{f}"]      \arrow[l, "{\text{id}}"]                   & X \arrow[u, hook]     
\end{tikzcd}
$$
At this point we introduce the line bundle on $\overline{X}'$ which we are going to work with.\ We define 
\begin{equation}\label{line bundle of concern}
F:=p_2^*\mathcal{L}\otimes (p_1^*\mathcal{L})^{\otimes(-q)}\in \text{Pic}(\overline{X}').
\end{equation}
Since the restriction $F|_{X_\eta}$ is isomorphic to the trivial line bundle by design,  \cite[ Corollaire 21.4.13]{grothendieck1964elements} implies that there is a line bundle $M$ on $S$ such that 
$$F|_{X}=\pi^*M.$$
We wish to extend  $M$ to $\overline{S}$.\ We work towards this goal in the next section.
\section{Extension of line bundle on the base}
In the previous section,\ we came upon the existence of a line bundle $M$ on $S$.\ We wish that $M$ comes from a line bundle $M'$ from above on $\overline{S}$.\ We can not wish to have this a priori, because we do not know that $\overline{S}$ is regular.\ However our idea here is to blow up $\overline{S}$ at certain closed subschemes lying outside $S$ so that we do not disturb $S$ and pass to a new $\overline{S}$ where this extension is possible.

We recall a basic extension result:
\begin{lemma}
\label{lemma:kunn}
Let us have a coherent sheaf $\mathcal{F}$ on an open sub-scheme $U$ of a noetherian scheme $X$.\ Also suppose $\mathcal{G}$ is a quasi-coherent sheaf on $X$ such that $\mathcal{F}\subseteq \mathcal{G}\mid_U$.\ Then there exists a coherent sheaf $\mathcal{F'}\subseteq \mathcal{G}$ on $X$ such that $\mathcal{F'}\mid_U=\mathcal{F}.$ 
\end{lemma}
\begin{proof}
See \cite[Exercise II.5.15]{hartshorne2013algebraic}.
\end{proof}
We apply the lemma to the following specific case:
\begin{lemma}
\label{lemma:effective}
If we have that $M$ is a coherent ideal sheaf on $S$,\ then we can find a coherent ideal sheaf $M'$ on $\overline{S}$ such that $M'\mid_S=M$.
\end{lemma}
\begin{proof}
We use Lemma \ref{lemma:kunn} with $M=\mathcal{F}$,\ $X=\overline{S}$,\ $U=S$ and $\mathcal{G}=\mathcal{O}_{\overline{S}}$.\ Note that the condition that $M$ is an ideal sheaf is necessary so that we can choose $\mathcal{G}$ to be the structure sheaf itself.
\end{proof}
It is trivial that we can always split a Weil divisor as a difference of two effective ones. The similar statement holds for Cartier divisors if the scheme is nice enough:  
\begin{lemma}
\label{lemma:liu}
Let $S$ be a quasi-projective scheme over a noetherian ring $k$ and $M$ be a line bundle on~$S$. Then we can find Cartier divisors $C, D$ and $E$ with $C, E$ being effective  such that
$\mathcal{O}_S(D)\cong M$ and $D$ is linearly equivalent to $C-E$.
\end{lemma}
\begin{proof}
Choose a Cartier divisor $D$ on $S$ such that $\mathcal{O}_S(D)\cong M$ (\cite[Proposition 7.1.32]{Liu}).\ By Lemma 7.1.31 of \cite{Liu}, there are very ample line bundles $L_1$ and $L_2$ such that $M\cong L_1\otimes L_2^{-1}$. By the very ampleness  there exist effective Cartier divisors $C$ and $E$ such that $\mathcal{O}_S(C)\cong L_1$ and $\mathcal{O}_S(E)\cong L_2$. Using \cite[Proposition 7.1.32]{Liu} again we get our claim.
\end{proof}
Now let us go back to the preparation of the proof of Theorem \ref{dynamical ST height inequality} and use the previous notations. Consider a Cartier divisor $D$ on $S$ such that $\mathcal{O}_S(D)\cong M$.\ Then by Lemma \ref{lemma:liu},\ we can choose to  write $D=C-E$ where $C$ and $E$ are effective divisors on $S$.\ Then we see from the construction of $\mathcal{O}_S(.)$ that both $\mathcal{O}_S(-C)$ and $\mathcal{O}_S(-E)$ are invertible ideal sheaves on $S$.\ Our idea is to extend these ideal sheaves with the help of Lemma \ref{lemma:effective} and blow up in their supports one by one.\ More formally,\ we have the following lemma and its proof explains the above procedure:
\begin{lemma}
\label{lemma:ext}
There is a projective variety  $S'$ and a birational morphism 
$$\gamma\colon S'\rightarrow \overline{S}$$
such that $\gamma\mid_{\gamma^{-1}(S)}$ is an isomorphism.\ Moreover if we identify $S$ as an open subscheme of $S'$ via $\gamma\mid_{\gamma^{-1}(S)}$ then we have a line bundle $M'$ on $S'$ such that $M'\mid_S\cong M$.
\end{lemma}
\begin{proof}
As explained in Lemma \ref{lemma:liu} above,\ using Lemma \ref{lemma:effective} we find a coherent ideal sheaf $\mathcal{C}$ on $\overline{S}$ such that $\mathcal{C}\mid_S=\mathcal{O}_S(-C)$.\  Now we consider the blow-up of $\overline{S}$
$$b_1\colon S_1\rightarrow \overline{S}$$
corresponding to the closed subscheme of the ideal sheaf $\mathcal{C}$.\ Note $b_1$ is an isomorphism when restricted to $S$ because the restriction of the ideal sheaf $\mathcal{C}$ to $S$ is $\mathcal{O}_S(-C)$ which is an invertible sheaf and blow-ups along invertible ideal sheaves are isomorphisms.\ Hence we can consider $S$ as an open subscheme of $S_1$ by identification via $b_1$.\ But now we have that the inverse image ideal sheaf $\overline{\mathcal{C}}=b_1^{-1}\mathcal{C}$ is an effective Cartier divisor by property of blow-ups (it is just the exceptional divisor corresponding to this blow up $b_1$).\ Note also that we have 
\begin{equation}\label{eq(1)}
\mathcal{O}_{S_1}(\overline{\mathcal{C}})\mid_S=\mathcal{O}_S(-C).
\end{equation}
Now via identification we view $S$ as an open sub-scheme of $S_1$ and $E$ is an effective Cartier divisor on $S$.\ Thus we once again have that $\mathcal{O}_S(-E)$ is an invertible ideal sheaf and we can once again apply Lemma \ref{lemma:effective} to obtain a coherent ideal sheaf extension $\mathcal{E}$ on $S_1$ such that $\mathcal{E}\mid_S=\mathcal{O}_S(-E)$.\ Using a similar argument as before we obtain that  we can safely blow-up in the associated closed sub-scheme without disturbing $S$.\ More formally we consider the blow-up 
$$b_2\colon S_2\rightarrow S_1$$along the closed subscheme corresponding to the coherent ideal sheaf $\mathcal{E}$.\ Then $b_2$ is an isomorphism over $S$ because the restriction of $\mathcal{E}$ to $S$ is the invertible ideal sheaf $\mathcal{O}_S(-E)$.\ At first note that as pullbacks of Cartier divisors are again Cartier divisors under dominant maps of integral schemes (so that is pullback is well defined),\ we get that $\tilde{\mathcal{C}}=b_2^*\overline{\mathcal{C}}$ is an invertible sheaf on $S_2$ with $\tilde{\mathcal{C}}\mid_S=\mathcal{O}_S(-C)$.\ Moreover,\ let us consider the effective Cartier divisor $\tilde{\mathcal{E}}$ which is the exceptional divisor corresponding to the blow up $b_2$.\ Then note that arguing like above and using identifications of $S$ via $b_2$,\ we get that 
\begin{equation}\label{eq(2)}
\mathcal{O}_{S_2}(\tilde{\mathcal{E}})\mid_S=\mathcal{O}_S(-E)
\end{equation}
 Now we are done because we take $M':=\mathcal{O}(\tilde{\mathcal{E}})\otimes\mathcal{O}(\tilde{\mathcal{C}})^{-1}\in \text{Pic}(S_2)$.\ Then from \eqref{eq(1)} and \eqref{eq(2)} we immediately see that 
$$M'\mid_S\cong M.$$ This completes the proof by taking $\gamma:=b_1\circ b_2$ and $S':=S_2$ (recall that blowing up preserves the projectivity).
\end{proof}
As we mentioned before,\ our goal is to pass to blow-ups where $M$ can be extended from $S$ to a suitable compactification of $S$.\ Using Lemma \ref{lemma:ext} we can find such blow-ups and without loss of generality we will replace the old $\overline{S}$ by $S_2$ in the lemma.\ Now we can assume that $M$ has an extension to $\overline{S}$.\ Note that this change of base space  does not effect the height on the base $S$ in Theorem \ref{dynamical ST height inequality} precisely because the blow-ups do not touch $S$\ (see Corollary \ref{corollary:change} in Section 5).
\begin{remark}
The above procedure could be easily done by an application of Hironaka's theorem on $\overline{S}$ and passing to a \emph{regular} compactification if the field $k$ is of characteristics zero.\ However, Hironaka's theorem is unknown for positive characteristic but blow-ups work just fine. This is one of the modifications we make from  \cite[Appendix A]{dimitrov2020uniformity}.
\end{remark}
\section{Height inequality using blow-ups}

In this section we begin by talking about the ideal sheaf of denominators of a Cartier divisor.\ Let us have a Cartier divisor given by the local equations $D=(U_i,f_i)$ where $U_i=\text{Spec}(A_i)$ are affine and $f_i\in \text{Frac}(A_i)^\times$ are the local equations of the Cartier divisor.\ Then we define the ideal sheaf of denominators $\mathcal{I}(D)$ to be 
$$\mathcal{I}(U_i)=\{b\in A_i\mid b\cdot f_i\in A_i\}$$
Note that clearly each $\mathcal{I}(U_i)$ is an ideal of $A_i$ and they glue together because $f_i$ and $f_j$ differ by an unit in the intersection $U_i\cap U_j$.\ Hence gluing these piece-wise sections we get a coherent ideal sheaf $\mathcal{I}$ which we call the ideal sheaf of denominators.

Now let us digress and prove a lemma which will be necessary for us later:
\begin{lemma}\label{lemma:nijo}
Let $X$ be an integral scheme,\ $U$ is an open subscheme of $X$ and $\mathcal{L}\in\emph{Pic}(X)$ such that $\mathcal{L}\mid_U$ is trivial.\ Then we can take a Cartier divisor $D$ on $X$ with $\mathcal{O}_X(D)\cong \mathcal{L}$ such that $D\mid_U$ is the trivial Cartier divisor on $U$ and hence $\mathcal{I}(D)\mid_U$ is trivial on $U$.
\end{lemma}
\begin{proof}
From the construction of $\mathcal{I}(D)$ it is clear to see that $\mathcal{I}(D)\mid_U=\mathcal{I}(D\mid_U)$.\ Hence to show the claim it is enough to choose a Cartier divisor corresponding to $\mathcal{L}$ such that its restriction to $U$ is trivial.\ We begin by choosing any Cartier representative $D$ of $\mathcal{L}$.\ Then as we have $\mathcal{L}\mid_U$ is trivial,\ we get that $D\mid_U$ is principal \ie there is an $f\in\kappa(X)$ such that $D\mid_U=(U,f)$.\ Then by replacing $D$ with $D'=D-(X,f)$ we get that $D'\mid_U=(U,1)$ which is the trivial Cartier divisor.\ However as we only change $D$ by a principal Cartier divisor,\ we stay in the same isomorphism class of $\mathcal{L}$.\ Hence we get that $\mathcal{L}\cong \mathcal{O}_X(D')$ and $\mathcal{I}(D')\mid_U$ along with $D'\mid_U$ is trivial as wanted. 
\end{proof}
Now we are ready to prove Theorem \ref{dynamical ST height inequality}. Recall the induced map $\overline{X}'\rightarrow \overline{S}$ which we also denote by $\pi$ for simplicity. Also recall
$$F\mid_{X}\cong \pi^*M$$ for some line bundle $M$ on $S$.\ This is where we will use the results of section 3.\ We can find a line bundle $M'$ on $\overline{S}$ such that $M'\mid_S\cong M$. In that case let us denote \begin{equation}\label{line bundle of concern2}
F':=F\otimes(\pi^*M')^{-1}.
\end{equation} In that case we can see that $F'\mid_{X}$ is trivial.\ Hence by Lemma \ref{lemma:nijo} above we can choose a Cartier divisor $D$ on $\overline{X}'$ such that $\mathcal{I}(D)\mid_{X}$ is trivial  and $\mathcal{O}_{\overline{X}'}(D)\cong F'$.\ Note also that by Lemma \ref{lemma:nijo},\ we can have $D\mid_{X}=(X,1)$,\ the trivial cartier divisor.\ Let us denote by $Z$ the closed subscheme defined by the ideal sheaf $\mathcal{I}(D)$.\ Then as $D$ restricted to $X$ is trivial,\ we get that $Z$ lies outside $X$.\ Let us blow up $\overline{X}'$ along $Z$.\ Note that as $Z$ lies outside $X$,\ this blow up 
$$\gamma_1\colon X_1\rightarrow \overline{X}'$$
is an isomorphism over $X$ and hence we are not changing our main object of interest (which was main point of the lemma above).\ Now as explained in the proof of \cite[Theorem 2.4]{fulton2013intersection} (Case 3), one has
$$\gamma_1^*D\simeq C-E$$ (linearly equivalence) for some effective Cartier divisor $C$ on $X_1$ and the exceptional divisor $E$ and we also have $C\cup E \subseteq \gamma_1^{-1}(\text{Supp}(D))$. 

For the sake of clarity,\ we re-iterate that we have the following diagram:
$$
\begin{tikzcd}
X_1 \arrow[r, "\gamma_1"] & \overline{X}' \arrow[d, "\pi"] \\
                                          & \overline{S}                            
\end{tikzcd}
$$
As $D\mid_{X}$ is trivial,\ we have that $\text{Supp}(D)\subseteq \overline{X}'-X$.\ Hence we get that the closed subset  $C\cup E$ lies outside $X$ if we identify $X$ with its inverse image via $\gamma_1$.\ Consider the scheme-theoretic union $C\cup' E$ in $X_1$ which gives a scheme structure on the closed subset $C\cup E$ (not necessarily reduced).\ Next consider the scheme-theoretic image of $C\cup' E$ under the morphism $\pi\circ \gamma_1$ and call it $\mathcal{Y}$.\ As the underlying set $C\cup E$ of the closed sub-scheme $C\cup' E$ is contained in $X_1-X$,\ we get that its scheme theoretic image $\mathcal{Y}$ is a closed sub-scheme whose underlying set is contained in $\overline{S}-S$.\\
Next we consider the scheme-theoretic inverse image ${\mathcal{Y}}_1:=(\pi\circ \gamma_1)^{-1}(\mathcal{Y})\subset X_1$. We consider the blow-up of $X_1$ along $\mathcal{Y}_1$
$$\gamma_2\colon X_2\rightarrow X_1.$$

As $\mathcal{Y}\cap S=\emptyset$,\ we have $\mathcal{Y}_1\cap X=\emptyset$.  Then note that the scheme-theoretic inverse image of $\mathcal{Y}$ under the dominant morphism $\tilde{\pi}:=\pi\circ(\gamma_1\circ \gamma_2)$ is also an effective Cartier divisor; let us call it $\tilde{\mathcal{Y}}$.\ Now we consider the change of the base space, \ie, we consider the blow up 
$$\gamma'\colon S'\rightarrow \overline{S}$$
along the closed subscheme $\mathcal{Y}$ of $\overline{S}$.\ As $\mathcal{Y}\subseteq \overline{S}-S$ we deduce that $\gamma'$ is an isomorphism over $S$.\  Then note that as the scheme-theoretic inverse image $\tilde{\mathcal{Y}}$ of $\mathcal{Y}$ under the morphism $\tilde{\pi}$ is an effective Cartier divisor,\ we get the dotted arrow making the following commutative  diagram by the universal property of blow-ups:
$$
\begin{tikzcd}
X_2 \arrow[r, "\gamma_2"] \arrow[rd, "\pi'", dotted] & X_1 \arrow[r, "\gamma_1"] & \overline{X}' \arrow[d, "\pi"] \\
                                                                & S' \arrow[r, "\gamma'"]                 & \overline{S}                           
\end{tikzcd}.
$$ 

Now we make a comment here that in the above diagram,\ each of the horizontal morphisms in the top row are isomorphisms over $X$.\ For $\gamma_1$ it is because the subscheme of denominators lies outside $X$ by Lemma \ref{lemma:nijo}  and for $\gamma_2$ it is because $\mathcal{Y}_1$ is contained in $X_1-X$, being the inverse image of a subset of $\overline{S}-S$.\ Hence we can view $X$ as an open dense subset of each of them.\ The morphism $\gamma'$ is also an isomorphism over $S$ as $\mathcal{Y}\subseteq\overline{S}-S$.

We use the same letter for an effective Cartier and its underlying closed subscheme  as we always do, but we will write $g^*$ to emphasize we are considering pullback of divisors ($g$ being a morphism of schemes such that we can define the pullback of a Cartier divisor). By the property of blow-ups,\ we get that $$\mathcal{Y}':=\gamma'^{-1}(\mathcal{Y})\subset S'$$ is an effective Cartier divisor.\ It is the exceptional divisor corresponding to the blow-up $\gamma'$.\  We claim that $$\pi'^*\mathcal{Y}'\pm\gamma_2^*(C-E)$$
is effective.\ Note it is enough to show that both 
 $$\gamma_2^*C\le\ \pi'^*(\mathcal{Y'})\quad  \text{and}\quad \gamma_2^*E\le \pi'^*(\mathcal{Y}').$$ 
By \cite[Corollary 11.49]{gortz2010algebraic}, the two conditions can be reformulated to that there are closed immersions for the underlying scheme structures
$$\gamma_2^{-1}(C)\xhookrightarrow{} \pi'^{-1}(\mathcal{Y}')\quad \text{and}\quad \gamma_2^{-1}(E)\xhookrightarrow{} \pi'^{-1}(\mathcal{Y}').$$ But these following conditions easily follow from diagram chasing and the construction.\ More precisely, there exists a closed immersion 
$$C\cup' E\xhookrightarrow{}(\pi\circ\gamma_1)^{-1}(\mathcal{Y})$$ where the right hand side is the scheme-theoretic inverse image.\ Hence we obtain that there exist closed immersions 
$$C\xhookrightarrow{}(\pi\circ\gamma_1)^{-1}(\mathcal{Y})\quad \text{and}\quad E\xhookrightarrow{}(\pi\circ\gamma_1)^{-1}(\mathcal{Y})$$ as $C$ and $E$ are closed subschemes of the scheme-theoretic union $C\cup' E$.\ As closed immersions are stable under base change we get that there are closed immersions 
$$\gamma_2^{-1}(C)\xhookrightarrow{}\tilde{\pi}^{-1}(\mathcal{Y})\quad \text{and}\quad \gamma_2^{-1}(E)\xhookrightarrow{}\tilde{\pi}^{-1}(\mathcal{Y}).$$ Now from the commutativity of the diagram above,\ we easily get that 
$\tilde{\pi}^{-1}(\mathcal{Y})=\pi'^{-1}(\mathcal{Y}')$ and the claim easily follows.

The effectiveness of the divisors $$\Lambda_{\pm}:=\pi'^*(\mathcal{Y}')\pm\gamma_2^*(C-E)$$ proves that outside their supports,\ the height is bounded by a constant.\ Rigorously speaking,\ this means that there is a constant $B>0$  that 
$$h_{\mathcal{Y}'}(\pi'(P))\pm h_{F'}(\gamma_1(\gamma_2(P)))\ge B$$
for all $P\in X_2(\overline{k})$ outside support of $\Lambda_{\pm}$.
Unwinding the definition \eqref{line bundle of concern2} of $F'$ and using the additivity and functoriality of the height machine,\ the above inequality can be rewritten as (up to bounded functions)
$$h_{\mathcal{Y}'}(\pi'(P))\pm h_{F}(\gamma_1(\gamma_2(P))\mp h_{M'}(\tilde{\pi}(P))\ge B.$$ This in particular means that 
\begin{equation}
    \label{ineq:1}
|h_{F}(\gamma_1(\gamma_2(P)))|\le |h_{\mathcal{Y}'}(\pi'(P))|+|h_{M'}(\tilde{\pi}(P))|+B.
\end{equation}

To conclude the inequality in Theorem \ref{dynamical ST height inequality}, it is sufficent to see that $X(\overline{k})\subset X_2(\overline{k})$ is away from the said basepoint-loci. More precisely, we claim that none of the points in $X(\overline{k})$ lies in the support of $\Lambda_{\pm}$.\   Observe that the supports of $\Lambda_{\pm}$ are given by $$\pi'^{-1}(\mathcal{Y'})\cup\gamma_2^{-1}(C\cup E)=\pi'^{-1}(\gamma'^{-1}(\mathcal{Y}))=\tilde{\pi}^{-1}(\mathcal{Y})\cup\gamma_2^{-1}(C\cup E)=\tilde{\pi}^{-1}(\mathcal{Y}).$$ 
Hence, if there is a point $P\in X(\overline{k})$ belongs to the support,\ then $\pi'(P)=\pi(P)\in \mathcal{Y}$.\ But  $X(\overline{k})$ has no intersection with $\gamma_2^{-1}(C\cup E)$. Thus no points of $X(\overline{k})$ can lie in support. Hence, after unwinding the definition of $F$, we deduce from \eqref{ineq:1} that (up to bounded functions) 
$$|h_{(f^*L\otimes L^{\otimes (-q)}}(P)|\le |h_{\mathcal{Y}}(\pi(P))|+|h_{M'}(\pi(P))|+B\ \text{for all}\ P\in X(\overline{k}).$$ The height machine  gives
\begin{equation}
\label{ineq:2}
    |h_{L}(f(P))-qh_{L}(P)|\le|h_{\mathcal{Y}}(\pi(P))|+|h_{M'}(\pi(P))|+B\ \text{for all}\ P\in X(\overline{k})
\end{equation} 
Note that inequality \eqref{ineq:2} is almost ready to be the inequality in Theorem \ref{dynamical ST height inequality} except that in the right hand side we have different terms $|h_{M'}|$ and $|h_{\mathcal{Y}}|$.

Each of the two heights can be bounded by some constant times the canonical Weil heights $h_{S'}$ and $h_{\overline{S}}$.\ The reason is that the standard Weil heights are given by very ample line bundles into their respective projective embeddings.\ This will also be made more clear in the next section.\  Hence to get the final height inequality we have to show that the height $h_{S'}$ is bounded by a constant times the height $h_{\overline{S}}$ when we restrict it to $X(\overline{k})$. 
\section{Last inequality}

This section serves as a bridging explanation for some assumptions we have made throughout the article.\ We have a recurring theme throughout the article where we pass via blow-up to a different compactification but the blow up does not touch certain open sub-schemes\ (that is an isomorphism over these open sub-schemes)\ which are typically $\mathcal{A}$ or $S$ for us.\ With help of Lemma \ref{lemma:z1} we show that the height functions only differ up to constant multiple factors when restricted to open sub-schemes over which these blow-ups are isomorphisms.\ On the light of the inequality in Theorem \ref{dynamical ST height inequality},\ this then justifies our change of base by application of blow-ups in section 3 and also shows us how to pass from inequality (2) in the previous section to the final inequality in Theorem \ref{dynamical ST height inequality}.

So let \[\pi: \mathbb{P}^{m'}\supset S'\to \overline{S}\subset \mathbb{P}^m\] be a blow up morphism, where $\overline{S}$ is a projective irreducible variety so that $S'$ is also projective and irreducible variety. Moreover, we have an open subset $S\subset \overline{S}$ which is not touched by the blow up, i.e.\ $\pi^{-1}(S)\cong S$. Consider two line bundles $L_1$ and $L_2$ on $S'$; let $h_1$, $h_2$ be the Weil heights attached to them respectively, while $h_{\overline{S}}$ (resp.\ $h_{S'}$) denotes the Weil height on $\overline{S}$ (resp.\ $S'$) induced by the standard height on $\mathbb{P}^m$ (resp.\ $\mathbb{P}^{m'}$) and the given closed immersion $\overline{S}\subset \mathbb{P}^m$ (resp.\ $S'\subset \mathbb{P}^{m'}$).  We want to bound $h_1$ and $h_2$ in terms of $h_{\overline{S}}$ (at least on a sub-domain).

For the next two assertions we adapt the convention that $h_n$ denotes the standard height on $\mathbb{P}^n,\, n\in \N$.
\begin{lemma}
\label{lemma:z1}
Let $V$ be a locally closed subvariety of $\mathbb{P}^n$ and suppose there is a morphism  $f: V\to \mathbb{P}^\ell$ defined over $k$ for some $\ell\in\mathbb{N}$. Then there exist constants $c_1, c_2>0$ such that
\[
\forall P\in V(\bar{k}): h_\ell(f(P))\leq c_1h_n(P)+c_2.
\]
\end{lemma}
\begin{proof}
See \cite[Proposition IV.1.7]{lang2013fundamentals}.
\end{proof}
\begin{corollary}
\label{lemma:z2}
Let $V$ be a variety over $k$. Let $c$ be an ample line bundle and $L$ any line bundle on $V$. Pick a height $h_c$ associated to $c$ such that $h_c\geq 0$ (possible because $c$ is ample). Then there exist constants $\gamma_1, \gamma_2>0$ such that
\[
\forall P\in V(\bar{k}): |h_L(P)|\leq \gamma_1h_c(P)+\gamma_2.
\]
\end{corollary}
\begin{proof}
(Cf.\ \cite[Proposition IV.5.4]{lang2013fundamentals}) By \cite[Lemma 7.1.31]{Liu} we can decompose $L=L_+\otimes L_-^{\otimes(-1)}$ where $L_+$ and $L_-$ are very ample line bundles on $V$ inducing closed immersions $V\hookrightarrow \mathbb{P}^{m_\pm}$ for some $m_\pm\in \mathbb{N}$. By the Height Machine we have \[h_L=h_{L_+}-h_{L-}=h_{m_+}-h_{m_-}\] as real-valued functions on $V(\bar{k})$ up to bounded functions. 

On the other hand, since we choose $h_c$ to be non-negative, by the Height Machine we may assume without loss of generality that $c$ is very ample, inducing a closed immersion $V\hookrightarrow \mathbb{P}^n$ for some $n\in \mathbb{N}$. Note that in this case $h_c=h_n$ as real-valued functions on $V(\bar{k})$ up to bounded functions.

Applying Lemma \ref{lemma:z1} to the immersions $V\hookrightarrow \mathbb{P}^{m_{\pm}}$, using the additivity of the Height Machine and rearranging the terms resulting from Lemma \ref{lemma:z1} (the right hand sides are both bounded by $\gamma h_c+O(1)$ for some $\gamma >0$), we are done.
\end{proof}
Let us come back to our beginning situation.\ We will take $L_1=\mathcal{O}_{S'}(\mathcal{Y})$ and $L_2=\gamma'^*M'$ which are line bundles on $S'$.(following our notations from previous sections)\  Applying Corollary \ref{lemma:z2} with $V=S'$, $c=\mathcal{O}_{\mathbb{P}^{m'}}(1)$, we are able to bound $h_1=h_{\mathcal{Y}}$ and $h_2=h_{\gamma'^*M'}$ with respect to $h_{S'}$.\ More precisely,\ using Corollary \ref{lemma:z2} we get positive constants $r_1,r_2,s_1,s_2$ such that 
\begin{equation}
    \label{ineq:3}
|h_{\mathcal{Y}}(P)|\le r_1\cdot h_{S'}(P)+s_1\ \text{and}\ |h_{M'}(P)|\le r_2\cdot h_{S'}(P)+s_2\ \text{for all}\ P\in S(\bar{k})
\end{equation}
where we have used the identification of $\gamma'^{-1}(S)$ and $S$ via $\gamma'$ and the functoriality of the height machine.

Now note that  $h_{S'}$ is non-negative. Since the algebraic points $P\in X(\overline{k})$ we concern are all mapped into $S(\overline{k})\subset \overline{S}(\overline{k})$  and $\gamma'$ is an isomorphism over $S$,\ the restrictions of $h_{S'}$ and $h_{\bar{S}}$ on $S(\bar{k})$ will differ only up-to constant multiples by Corollary \ref{corollary:change} below.\ Hence by inequality (3) above we can absorb all the terms in inequality (2) of last section and arrive at Theorem \ref{dynamical ST height inequality}.
\begin{corollary}
\label{corollary:change}
Let us have  two projective varieties $S'$ and $\overline{S}$ embedded in $\mathbb{P}_k^m$ and $\mathbb{P}_k^n$ respectively and let $h_m$ and $h_n$ be the standard Weil heights of these projective spaces respectively.\ Furthermore,\ suppose we have an open sub-scheme $S\subseteq \overline{S}$ and a morphism 
$$\pi\colon S'\rightarrow \overline{S}$$
such that $\pi\mid_{\pi^{-1}(S)}$ is an isomorphism onto $S$.\ Then we have positive constants $c_1,c_2,d_1,d_2$ such that 
$$c_1\cdot h_n(P)-d_1\le h_m(P)\le c_2\cdot h_n(P)+d_2$$
for all $P\in S(\overline{k})$ (where we identify $\pi^{-1}(S)$ and $S$ via $\pi$)
\end{corollary}
\begin{proof}
The proof consists of an application of Lemma \ref{lemma:z1} once with $V=\pi^{-1}(S)$ and $$f=\pi\mid_{\pi^{-1}(S)}\colon \pi^{-1}(S)\rightarrow \mathbb{P}_k^n$$ and once in the reverse direction with $V=S$ and $$f=(\pi\mid_{\pi^{-1}(S)})^{-1}\colon S\rightarrow \mathbb{P}_k^m$$ Note that for the application of the second case,\ it is vital that $\pi$ is an isomorphism over $S$.\ Once we write down the two inequalities obtained this way and after some algebraic manipulations.\ we clearly get the inequalities of the Lemma.
\end{proof}
\begin{proof}[Proof of Theorem \ref{dynamical ST height inequality}]
Note that this lemma essentially shows that if we pass to a blow up of $\overline{S}$ which does not touch $S$,\ the height functions (coming from the projective embeddings of $\overline{S}$ and $S'$) are equivalent to each other up to some universal constant multiples.\ In the view of the role of the height on the base space in the equality of Theorem \ref{dynamical ST height inequality},\ this clearly justifies the change of base by blow-up that we did in Section 3 because in the notation of the previous sections,\ we have $h_m=h_{S'}$ and $h_n=h_{\overline{S}}$.

Now we come to obtaining final inequality in Theorem \ref{dynamical ST height inequality} from inequality (\ref{ineq:2}).\ Using inequality (\ref{ineq:3}) in inequality (\ref{ineq:2}),\ we get (change the constant $B$ from the last section if necessary)
$$|h_{L}(f(P))-qh_{L}(P)|\le(r_1+r_2)\cdot h_{S'}(\pi(P))+(s_1+s_2)+B\  \text{for all}\ P\in X(\overline{k})$$
Next we apply the inequality of Corollary \ref{corollary:change} noting that $h_m=h_{S'}$ and $h_n=h_{\overline{S}}$ to get that 
$$|h_{L}(f(P))-qh_{L}(P)|\le c_2(r_1+r_2)h_{\overline{S}}(\pi(P))+d_2(r_1+r_2)+(s_1+s_2)\  \text{for all}\ P\in X(\overline{k})$$
Now finally from the above we have 
$$|h_{L}(f(P))-qh_{L}(P)|\le c\cdot \max\{1,h_{\overline{S}}(\pi(P))\}\ \text{for all}\ P\in X(\overline{k})$$
where we take $c=2\cdot \text{max}\{c_2(r_1+r_2), d_2(r_1+r_2)+(s_1+s_2)\}$ which finally proves Theorem \ref{dynamical ST height inequality} as $c$ does not depend on $P$ clearly.
\end{proof}
\begin{remark}
Both Theorem \ref{dynamical ST height inequality} and Corollary \ref{ST's height ineq} are natural generalization of standard facts of the theory of heights on a variety (over a point) to the relative setting. 

More precisely, if $S=\text{Spec}(k)$ is a single point and $X$ is an abelian variety over $S$, then Corollary \ref{ST's height ineq} reflects the basic fact that the N\'eron-Tate height $\hat{h}_X$ and the Weil height $h_{\overline{X}}$ lie in the same class of the height function associated to the line bundle $L$ (cf. Weil's height machine) that gives rise to the projective embedding of $X$. 

Still consider an abelian variety $X$ over a point $S$ and let $L$ be the (even and amplie) line bundle~$L$ as before. It is also well-known that the difference $|h_L([2]P)-4h_L(P)|$ is bounded and for the induced N\'eron-Tate height $\hat{h}$ one has a genuine equality $\hat{h}([2]P)=4\hat{h}(P)$  (for all $P\in X(\overline{k})$). Now in the relative (dynamical) setting, if $(X,f,L,q)_{/S}$ is a polarized dynamical system, where $S$ is a quasi-projective $k$-variety, then one can construct a so-called $f$-\emph{canonical height} $h_f$ such that for all algebraic points $P\in X(\overline{k})$, one has 
\[
h_f(f(P))=qh_f(P).
\] See \cite[Proposition 6.1.4(2)]{YZ21}, which is a generalization of \cite[Theorem 2.4(a)]{zhang1995small}.
\end{remark}
\section{Comparison of generic and special N\'eron-Tate Heights}
In this section we consider the problem of estimating exceptional points in the specialization map of a family of abelian varieties. Suppose we have a projective morphism of varieties $\pi\colon\overline{\mathcal{A}}\rightarrow \mathbb{P}_k^m$ \emph{i.e} there exists a closed immersion $\overline{\mathcal{A}}\xhookrightarrow{i}\mathbb{P}_k^n\times \mathbb{P}_k^m$ such that composing it with the projection to $\mathbb{P}_k^m$ gives $\pi$. We further assume that the generic fiber $\overline{\mathcal{A}}_{\eta}\rightarrow k(\mathbb{P}_m)$ is an abelian variety over the function field $k(\mathbb{P}_k^m)$. Note that here we do not assume that the base space $\mathbb{P}_k^m$ is a curve as in Theorem B of \cite{Sil83} neither will we assume that $\pi$ is flat as in Theorem A in \cite{Sil83}. Furthermore we will assume the base field $k$ to be a global field of any characteristic. Since the generic fiber $\overline{\mathcal{A}}_{\eta}$ is an abelian variety, using standard results we can deduce that there exists an open dense subset $S$ of $\mathbb{P}_k^m$ such that $\pi\colon \pi^{-1}(S)\rightarrow S$ is an abelian scheme and we set $A:=\pi^{-1}(S)$. Since $S$ already comes prescribed with an immersion inside $\mathbb{P}_k^m$, we also denote $\overline{S}:=\mathbb{P}_k^m$ as a projective compactification. Furthermore recall that we have the very ample line bundle $\mathcal{L}=\mathcal{O}(1,1)|_{\overline{\mathcal{A}}}$ and we assume that $L:=\mathcal{L}|_{\overline{A}}$ is an even relative ample line bundle on the abelian scheme $\mathcal{A}$. Then this constitutes a polarised dynamical system as described in Example \ref{abelian scheme as dynamical system} and subsequently we have certain height functions. We denote the canonical height on $\overline{\mathcal{A}}_{\eta}(K)$ corresponding to the even ample line bundle $\mathcal{L}_{\eta}$ as $\hat{h_{\mathcal{L}_{\eta}}}$ and similarly the canonical fiber-wise  height on $\mathcal{A}(k)$ corresponding to the even ample line bundle $L:= \mathcal{L}|_{\overline{A}}$ as $\hat{h_{\overline{\mathcal{L}}}}$. Moreover we denote the naive Weil height corresponding to $\mathcal{L}$ coming from the immersion $\overline{\mathcal{A}}\xhookrightarrow{} \mathbb{P}_k^n\times \mathbb{P}_k^m$ as $h_{\mathcal{L}}$ and the usual Weil height of the projective space $\mathbb{P}_k^m$ by $h_m$.\\
Now suppose we have a closed rational point $P\in\overline{\mathcal{A}}_{\eta}$ \emph{i.e} a closed immersion 
$\text{Spec}(k(S))\xhookrightarrow{P}\overline{\mathcal{A}}_{\eta}\rightarrow \text{Spec}(k(S))$ which is a section. Then since $S$ is normal, Weil's extension theorem\break (\cite[Sec 4.4, Thm 1]{Bosch} together with Hartog's extension theorem implies that $P$ extends to a section $S\xhookrightarrow{}\mathcal{A}\rightarrow S$ which we denote by $P_S$. Hence for any closed point $t\in S(k)$, we can define the specialisation map $P\to P_t$ which is just the base change as in the cartesian diagram 
\[\begin{tikzcd}
\label{figure:one}
\text{Spec}(k(t)) \arrow[r, "P_t", hook, shift right] \arrow[d] & \mathcal{A}_t \arrow[d] \\
S \arrow[r, "P_S", hook]                                               & \mathcal{A}            
\end{tikzcd}\]
where $k(t)$ denotes the field of definition of $t$ over $k$
and we are interested in all those $t$ such that this sepcialisation map is not injective. We begin with a limit theorem analogous to \cite[Theorem 4.1]{SilCall}.
\begin{theorem}
\label{theorem:limittheorem}
Suppose $\overline{A}\rightarrow \mathbb{P}_k^m$ be a projective morphism such that the generic fiber $\overline{\mathcal{A}}_{\eta}$ over the function field $K:=k(\mathbb{P}_k^m)$ is an abelian variety. Furthermore suppose $S$ is the maximal open dense subset of $\mathbb{P}_k^m$ such that $\pi\colon\pi^{-1}(S)\rightarrow S$ is an abelian scheme and suppose $H:=\mathbb{P}_k^m\backslash S$. Then for all Zariski closed sub-varieties of $\mathbb{P}_k^m$ with $Z\cap H=\phi$, we have
\[\lim_{h_m(t)\to \infty, t\in Z(k(t))}\frac{\hat{h_{\mathcal{L}}}(P_t)}{h_m(t)}=\hat{h_{\mathcal{L}_{\eta}}}(P)\]

\end{theorem}
\begin{proof}
Consider the \say{spreaded out} section appearing in the lower horizontal arrow of the cartesian square in Figure \ref{figure:one}. We begin by noting that $P_S\colon S\rightarrow \mathcal{A}\subset \overline{\mathcal{A}}$ induces a rational morphism from $\mathbb{P}_k^m$ to $\mathbb{P}_k^n$ via the sequence 
$$
\begin{tikzcd}
S\subseteq \mathbb{P}_k^m \arrow[r, "P_S", dotted, hook] \arrow[rrd, "f_P", dotted] & \mathcal{A} \arrow[r, hook] & \mathbb{P}_k^m\times \mathbb{P}_k^n \arrow[d] \\
                                                                                    &                             & \mathbb{P}_k^n                               
\end{tikzcd}
$$
which is regular over $S$. If we denote by $h_{\mathcal{L}}$ the Weil height of $\mathcal{L}$, then we have 
\[h_{\mathcal{L}}(P_t)=h_m(P_S(t))+h_m(t)+O(1)\]
where $h_m$ denotes the Standard Weil height in $P_k^m$.  Note that from Theorem \ref{dynamical ST height inequality} we can conclude that 
\[\hat{h_{\mathcal{L}}}(P_t)=h_{\mathcal{L}}(P_t)+O(h_m(t))\ \text{and}\ \hat{h_{\mathcal{L}_{\eta}}}(P)=h_{\mathcal{L}_{\eta}}(P)+O(1)\ \text{as}\ h_m(t)\to\infty\ \text{and}\ t\in S(k)\]
where $h_{\mathcal{L}_{\eta}}$ is the Weil height of $\mathcal{L}_{\eta}$.
as $h_m(t)\to\infty$. Suppose we manage to show 
\begin{equation}
\label{eqn:ratio}
\frac{h_m(f_P(t))}{h_m(t)}=h_{\mathcal{L}_{\eta}}(P)+o(1)\ \text{as}\ h_m(t)\to\infty\ \text{and}\ t\in Z(k)
\end{equation}
Then clearly from the above three inequalities, we can find a constant independent of $t\in Z(k)$ such that 
\[|\lim_{h_m(t)\to \infty,\ t\in Z(k)}\frac{\hat{h_{\mathcal{L}}}(P_t)}{h_m(t)}-\hat{h_{\mathcal{L}_{\eta}}}(P)|\le C\]
for all $t\in Z(k)$. Noting the homogeneity of both the terms in L.H.S of the above with the morphism $[N]$ clearly shows the claim. Hence we are reduced to showing \ref{eqn:ratio}.\\
To this end suppose that the rational map $f_P\colon \mathbb{P}_k^m\dashrightarrow \mathbb{P}_k^n$ is given by homogenous polynomials of degree $d$ in the projective co-ordinates of $\mathbb{P}_k^m$. Then by the discussion at the end of section 2.4 of \cite{bombieri2007heights} we deduce that $h_{\mathcal{L}_{\eta}}(P)=d$ as $X=\mathbb{P}_k^m$. On the other hand from Lemma 1.6 of \cite{lang2013fundamentals} we have
\[h_m(f_P(t))=d\cdot h_m(t)+O(1)=h_{\mathcal{L}_{\eta}}(P)\cdot h_m(t)+O(1)\ \text{as}\ h_m(t)\to\infty\ \text{and}\ t\in Z(k)\]
as the domain $U=\mathbb{P}_k^m\backslash H$ of $f_P$ properly contains $Z$ which clearly shows \ref{eqn:ratio} as $h_m(t)\to \infty$ and finishes our proof.
\end{proof}
 Note that for all $t\in S(k)$, we obtain a map 
\[\sigma_t\colon \mathcal{A}_{\eta}(K)\rightarrow \mathcal{A}_t(k(t))\]
\[P\mapsto P_t\]
Then with the above Theorem, arguing exactly similarly as in the proof of Theorem C in \cite{Sil83} we obtain a variant of the same.
\begin{corollary}
\label{coroll:injectivity}
Suppose $\pi\colon\overline{\mathcal{A}}\rightarrow \mathbb{P}_k^m$ be as in Theorem \ref{theorem:limittheorem} and furthermore suppose that the $K/k$-trace $\emph{tr}_{K/k}(\mathcal{A_{\eta}})=0$. Then there exists a Zariski closed subset $H\subseteq \mathbb{P}_k^m$ such that for all Zariski closed subvarieties $Z\subseteq \mathbb{P}_k^m$ with $Z\cap H=\phi $ the set 
\[\{t\in Z(k)\mid \sigma_t\ \text{is not injective}\}\]
is of bounded height with respect to $h_m(\cdot)$.
\end{corollary}
\begin{remark}
The results of this chapter extend those of \cite{Sil83} and \cite{SilCall} to the case where the parameter space is of higher dimension. But we only consider the case when the family is parametrised by a projective space and not any arbitary projective variety. Furthermore we restrict our attention on closed sub-varieties of $\mathbb{P}_k^m$ contained in $S$ rather than all of $S$ itself as done in \cite{Sil83}. A full generalisation of Theorem \ref{theorem:limittheorem} to higher dimensions is not possible as can be seen by easy counter examples. Note that Theorem \ref{theorem:limittheorem} is true in the dynamical setup, as our height inequality (Theorem \ref{dynamical ST height inequality}) is.  We leave the obvious modifications to the reader.
\end{remark}
\printbibliography

@article {grothendieck1964elements,
  title = {Éléments de géométrie algébrique. IV. Étude locale des schémas et des
morphismes de schémas I-IV},
  author = {Grothendieck, A.},
  journal = {Inst. Hautes Études Sci. Publ. Math.},
  volume = {20,24,28,32},
  number = {},
  pages = {},
  year = {1964-1967},
  publisher = {Springer}
}

@article{dimitrov2020uniformity,
     AUTHOR = {Dimitrov, V. and Gao, Z. and Habegger, P.},
     TITLE = {Uniformity in {M}ordell-{L}ang for curves},
   JOURNAL = {Ann. of Math. (2)},
  FJOURNAL = {Annals of Mathematics. Second Series},
    VOLUME = {194},
      YEAR = {2021},
    NUMBER = {1},
     PAGES = {237--298},
      ISSN = {0003-486X},
   MRCLASS = {11G50 (11G30 14G05 14G25)},
  MRNUMBER = {4276287},
       DOI = {10.4007/annals.2021.194.1.4},
       URL = {https://doi.org/10.4007/annals.2021.194.1.4},
}

@book{gortz2010algebraic,
    AUTHOR = {G\"{o}rtz, U. and Wedhorn, T.},
     TITLE = {Algebraic geometry {I}},
    SERIES = {Advanced Lectures in Mathematics},
      NOTE = {Schemes with examples and exercises},
 PUBLISHER = {Vieweg + Teubner, Wiesbaden},
      YEAR = {2010},
     PAGES = {viii+615},
      ISBN = {978-3-8348-0676-5},
   MRCLASS = {14-01},
  MRNUMBER = {2675155},
MRREVIEWER = {C\'{\i}cero Carvalho},
       DOI = {10.1007/978-3-8348-9722-0},
       URL = {https://doi.org/10.1007/978-3-8348-9722-0},
}

@book{fulton2013intersection,
    AUTHOR = {Fulton, W.},
     TITLE = {Intersection Theory},
    SERIES = {Ergebnisse der Mathematik und ihrer Grenzgebiete. 3. Folge. A
              Series of Modern Surveys in Mathematics},
    VOLUME = {2},
   EDITION = {Second},
 PUBLISHER = {Springer-Verlag, Berlin},
      YEAR = {1998},
     PAGES = {xiv+470},
      ISBN = {3-540-62046-X; 0-387-98549-2},
   MRCLASS = {14C17 (14-02)},
  MRNUMBER = {1644323},
       DOI = {10.1007/978-1-4612-1700-8},
       URL = {https://doi.org/10.1007/978-1-4612-1700-8},
}

@book{lang2013fundamentals,
    AUTHOR = {Lang, S.},
     TITLE = {Fundamentals of {D}iophantine geometry},
 PUBLISHER = {Springer-Verlag, New York},
      YEAR = {1983},
     PAGES = {xviii+370},
      ISBN = {0-387-90837-4},
   MRCLASS = {11-02 (11Dxx 11Gxx 14G25)},
  MRNUMBER = {715605},
MRREVIEWER = {Gerd Faltings},
       DOI = {10.1007/978-1-4757-1810-2},
       URL = {https://doi.org/10.1007/978-1-4757-1810-2},
       note={2nd Edition},
}

@book{Liu,
    AUTHOR = {Liu, Q.},
     TITLE = {Algebraic geometry and arithmetic curves},
    SERIES = {Oxford Graduate Texts in Mathematics},
    VOLUME = {6},
      NOTE = {Translated from the French by Reinie Ern\'{e},
              Oxford Science Publications},
 PUBLISHER = {Oxford University Press, Oxford},
      YEAR = {2006},
     PAGES = {xvi+576},
      ISBN = {0-19-850284-2},
   MRCLASS = {14-01 (11G30 14A05 14A15 14Gxx 14Hxx)},
  MRNUMBER = {1917232},
MRREVIEWER = {C\'{\i}cero Carvalho},
}

@book{hartshorne2013algebraic,
AUTHOR = {Hartshorne, R.},
TITLE = {Algebraic Geometry},
PUBLISHER = {Springer},
    SERIES = {Graduate Texts in Mathematics},
    volume={52},
YEAR ={1977},
note ={Corrected 14th printing} 
}

@book{bombieri2007heights,
    AUTHOR = {Bombieri, E. and Gubler, W.},
     TITLE = {Heights in {D}iophantine geometry},
    SERIES = {New Mathematical Monographs},
    VOLUME = {4},
 PUBLISHER = {Cambridge University Press, Cambridge},
      YEAR = {2006},
     PAGES = {xvi+652},
      ISBN = {978-0-521-84615-8; 0-521-84615-3},
   MRCLASS = {11G50 (11-02 11G10 11G30 11J68 14G40)},
  MRNUMBER = {2216774},
MRREVIEWER = {Yuri Bilu},
       DOI = {10.1017/CBO9780511542879},
       URL = {https://doi.org/10.1017/CBO9780511542879},
}

@article {Sil83,
    AUTHOR = {Silverman, J. H.},
     TITLE = {Heights and the specialization map for families of abelian
              varieties},
   JOURNAL = {J. Reine Angew. Math.},
  FJOURNAL = {Journal f\"{u}r die Reine und Angewandte Mathematik. [Crelle's
              Journal]},
    VOLUME = {342},
      YEAR = {1983},
     PAGES = {197--211},
      ISSN = {0075-4102},
   MRCLASS = {14K15 (14D10 14G25)},
  MRNUMBER = {703488},
MRREVIEWER = {Gerd Faltings},
       DOI = {10.1515/crll.1983.342.197},
       URL = {https://doi.org/10.1515/crll.1983.342.197},
}

@article{gao2019heights,
      title={Heights in families of abelian varieties and the Geometric Bogomolov Conjecture}, 
      author={Gao, Z. and Habegger, P.},
     TITLE = {Heights in families of abelian varieties and the geometric
              {B}ogomolov conjecture},
   JOURNAL = {Ann. of Math. (2)},
  FJOURNAL = {Annals of Mathematics. Second Series},
    VOLUME = {189},
      YEAR = {2019},
    NUMBER = {2},
     PAGES = {527--604},
      ISSN = {0003-486X},
   MRCLASS = {11G10 (11G50 14G25 14G40 14K15)},
  MRNUMBER = {3922127},
MRREVIEWER = {John L. Boxall},
       DOI = {10.4007/annals.2019.189.2.3},
       URL = {https://doi-org.ezproxy.leidenuniv.nl/10.4007/annals.2019.189.2.3},
}

@book {Laz04,
    AUTHOR = {Lazarsfeld, R.},
     TITLE = {Positivity in algebraic geometry. {I}},
    SERIES = {Ergebnisse der Mathematik und ihrer Grenzgebiete. 3. Folge. A
              Series of Modern Surveys in Mathematics [Results in
              Mathematics and Related Areas. 3rd Series. A Series of Modern
              Surveys in Mathematics]},
    VOLUME = {48},
      NOTE = {Classical setting: line bundles and linear series},
 PUBLISHER = {Springer-Verlag, Berlin},
      YEAR = {2004},
     PAGES = {xviii+387},
      ISBN = {3-540-22533-1},
   MRCLASS = {14-02 (14C20)},
  MRNUMBER = {2095471},
MRREVIEWER = {Mihnea Popa},
       DOI = {10.1007/978-3-642-18808-4},
       URL = {https://doi.org/10.1007/978-3-642-18808-4},
}

@book{YZ21,
      title={Adelic line bundles over quasi-projective varieties}, 
      author={Yuan, X. and Zhang, S.-W.},
      year={2022},
      note={Available at \url{https://arxiv.org/abs/2105.13587}},
      archivePrefix={arXiv},
      primaryClass={math.NT},
}

@article {Holmes,
    AUTHOR = {Holmes, D.},
     TITLE = {Torsion points and height jumping in higher-dimensional
              families of abelian varieties},
   JOURNAL = {Int. J. Number Theory},
  FJOURNAL = {International Journal of Number Theory},
    VOLUME = {15},
      YEAR = {2019},
    NUMBER = {9},
     PAGES = {1801--1826},
      ISSN = {1793-0421},
   MRCLASS = {14G40 (11G50)},
  MRNUMBER = {4015515},
MRREVIEWER = {Nathan Kaplan},
       DOI = {10.1142/S179304211950101X},
       URL = {https://doi-org.ezproxy.leidenuniv.nl/10.1142/S179304211950101X},
}

@article{zhang1995small,
  title={Small points and adelic metrics},
  author={Zhang, S.-W.},
  journal={Journal of Algebraic Geometry},
  volume={4},
  number={2},
  pages={281--300},
  year={1995},
  publisher={Citeseer}
}

@article {SilCall,
    AUTHOR = {Call, Gregory S. and Silverman, Joseph H.},
     TITLE = {Canonical heights on varieties with morphisms},
   JOURNAL = {Compositio Math.},
  FJOURNAL = {Compositio Mathematica},
    VOLUME = {89},
      YEAR = {1993},
    NUMBER = {2},
     PAGES = {163--205},
      ISSN = {0010-437X},
   MRCLASS = {11G35 (14G25)},
  MRNUMBER = {1255693},
MRREVIEWER = {Takeshi Ooe},
       URL = {http://www.numdam.org/item?id=CM_1993__89_2_163_0},
}

@book {Bosch,
    AUTHOR = {Bosch, Siegfried and L\"{u}tkebohmert, Werner and Raynaud, Michel},
     TITLE = {N\'{e}ron models},
    SERIES = {Ergebnisse der Mathematik und ihrer Grenzgebiete (3) [Results
              in Mathematics and Related Areas (3)]},
    VOLUME = {21},
 PUBLISHER = {Springer-Verlag, Berlin},
      YEAR = {1990},
     PAGES = {x+325},
      ISBN = {3-540-50587-3},
   MRCLASS = {14K15 (11G10 14L15)},
  MRNUMBER = {1045822},
MRREVIEWER = {James Milne},
       DOI = {10.1007/978-3-642-51438-8},
       URL = {https://doi.org/10.1007/978-3-642-51438-8},
}
\end{document}